\documentclass[12pt]{article}

\usepackage        {anysize }%
\usepackage        {latexsym}%
\usepackage        {amsmath }%
\usepackage        {amssymb }% 
\usepackage        {amsthm  }%
\usepackage        {paralist}%
\usepackage        {url     }% 

\DeclareMathOperator{\bier}{Bier}%
\DeclareMathOperator{\rank}{rk  }%
\DeclareMathOperator{\stellar}{sd}%

\newcommand  {\R  }{\mathbb R  }%
\newcommand  {\ov }{\overline  }%
\newcommand  {\sm }{{\setminus}}%

\theoremstyle{plain}%
\newtheorem{theorem}         {Theorem}    [section]%
\newtheorem{cor}    [theorem]{Corollary}%
\newtheorem{lemma}  [theorem]{Lemma}%
\newtheorem{prop}   [theorem]{Proposition}%
\newtheorem{defi}   [theorem]{Definition}%
\theoremstyle{remark}
\newtheorem*{remark*}{Remark}%

\reversemarginpar

\begin{document}
\title{{\huge\bfseries Bier spheres and posets}}
\author{{\sc Anders Bj\"orner}\thanks{%
        Research partially supported by the European Commission's IHRP
    Programme, grant HPRN-CT-2001-00272, ``Algebraic Combinatorics in Europe''}\\
        \small Dept.\ Mathematics\\[-1.7mm]
        \small KTH Stockholm\\[-1.7mm]
        \small S-10044 Stockholm, Sweden\\[-1.7mm]
        \small \url{bjorner@math.kth.se}
\and    {\sc Andreas Paffenholz}\setcounter{footnote}{6}\thanks{%
    Research supported by the Deutsche Forschungsgemeinschaft
    within the European graduate program ``Combinatorics, Geometry, and
    Computation'' (No. GRK 588/2)}\\
        \small Inst.\ Mathematics, MA 6-2\\[-1.7mm]
        \small TU Berlin\\[-1.7mm]
        \small D-10623 Berlin, Germany\\[-1.7mm]
        \small \url{paffenholz@math.tu-berlin.de}
\and   {\sc Jonas Sj\"ostrand        }\\
        \small Dept.\ Mathematics\\[-1.7mm]
        \small KTH Stockholm\\[-1.7mm]
        \small S-10044 Stockholm, Sweden\\[-1.7mm]
        \small \url{jonass@math.kth.se}
\and{\sc G\"unter M. Ziegler}$^*$\setcounter{footnote}{6}\thanks{$^*$Partially 
          supported by Deutsche Forschungsgemeinschaft, via the
  DFG Research Center ``Mathematics in the Key Technologies'' (FZT86),
  the Research Group ``Algorithms, Structure, Randomness'' (Project ZI 475/3),
  a Leibniz grant (ZI 475/4),
     and by the German Israeli Foundation (G.I.F.)}\\
        \small Inst.\ Mathematics, MA 6-2\\[-1.7mm]
        \small TU Berlin\\[-1.7mm]
        \small D-10623 Berlin, Germany\\[-1.7mm]
        \small \url{ziegler@math.tu-berlin.de}
        }
\date{~\\[3mm]{\small April 9, 2004}\\[3mm]
 Dedicated to Louis J. Billera on occasion of his 60th birthday}
\maketitle

\begin{abstract}\noindent%
  In 1992 Thomas Bier presented a strikingly simple method to produce a
  huge number of simplicial $(n-2)$-spheres on $2n$ vertices as deleted
  joins of a simplicial complex on $n$ vertices with its combinatorial
  Alexander dual.
  
  Here we interpret his construction as giving the poset of all the
  intervals in a boolean algebra that ``cut across an ideal.''  Thus we
  arrive at a substantial generalization of Bier's construction: the
  \emph{Bier posets} $\bier(P,I)$ of an arbitrary bounded poset $P$ of
  finite length. In the case of face posets of PL spheres this yields
  cellular ``generalized Bier spheres.''  In the case of Eulerian or
  Cohen-Macaulay posets $P$ we show that the Bier posets $\bier(P,I)$
  inherit these properties.
  
  In the boolean case originally considered by Bier, we show that all the
  spheres produced by his construction are shellable, which yields ``many
  shellable spheres,'' most of which lack convex realization.
  Finally, we present simple explicit formulas for the $g$-vectors of
  these simplicial spheres and verify that they satisfy a strong form of
  the $g$-conjecture for spheres.
\end{abstract}

%%%%%%%%%%%%%%%%%%%%%%%%%%%%%%%%%%%%%%%%%%%%%%%%%%%%%%%%%%%%%%%%%%%%%

\section*{Introduction}

In unpublished notes from 1992, Thomas Bier \cite{Bier92} described a
strikingly simple construction of a large number of simplicial PL spheres.
His construction associates a simplicial $(n-2)$-sphere with
$2n$ vertices to any simplicial complex $\Delta\subset2^{[1,n]}$ on 
$n$ vertices (here $[1,n]:=\{1,2,\dots , n\}$), 
by forming the ``deleted join'' of the complex $\Delta$ with its
combinatorial Alexander dual, $\Delta^*:=\{F\subset[1,n]:[1,n]\sm
F\notin\Delta\}$.  Bier proved that this does indeed yield PL spheres by
verifying that any addition of a new face to $\Delta$ amounts to a
bistellar flip on the deleted join of $\Delta$ with its Alexander dual
$\Delta^*$.  A short published account of this proof is given in
Matou\v{s}ek \cite[Sect.~5.6]{Matousek03}, to where we also
refer for the definition of deleted joins.
See de Longueville
\cite{longueville:_bier} for a simple alternative proof.  
\smallskip

\noindent
In this paper we generalize and further analyze Bier's construction:
\begin{itemize}\itemsep=0pt
\item We define 
more general ``Bier posets'' $\bier(P,I)$, where $P$
  is an arbitrary bounded poset of finite length and $I\subset P$ is an
  order ideal.
\item We show that in this generality, the order complex of $\bier(P,I)$
  is PL homeomorphic to that of $P$: It may be obtained by a sequence of
  stellar subdivisions of edges.
\item If $P$ is an Eulerian or Cohen-Macaulay poset or lattice, then
  $\bier(P,I)$ will have that property as well.
\item If $P$ is the face lattice of a regular CW PL-sphere $\mathcal S$,
  then the lattices $\bier(P,I)$ are again face lattices of regular CW
  PL-spheres, the ``Bier spheres'' of $\mathcal S$.
\item In the case of the $(n-1)$-simplex, where $P=B_n$ is a boolean
  algebra, and the ideal in $B_n$ may be interpreted as an abstract
  simplicial complex $\Delta$, one obtains the ``original'' Bier spheres
  as described in \cite{Bier92}, with face lattice $\bier(B_n,\Delta)$.
  We prove that all these simplicial spheres are shellable.
\item The number of these spheres is so great that for large $n$ most of
  the Bier spheres $\bier(B_n,\Delta)$ are not realizable as polytopes.
  Thus Bier's construction provides ``many shellable spheres'' in the
  sense of Kalai \cite{Kalai88} and Lee \cite{Lee00}; see
  also \cite[p.~116]{Matousek03}.
  Similarly, for special choices of the simplicial complex
  $\Delta$ in $B_n$, and even $n$, we obtain many nearly neighborly centrally
  symmetric $(n-2)$-spheres on $2n$ vertices.
\item The $g$-vector of a Bier sphere 
  $\bier(B_n,\Delta)$ can be expressed explicitly in
  terms of the $f$-vector of $\Delta$.  We show that
these $g$-vectors actually are $K$-sequences, and thus 
  they satisfy a strong form of the $g$-conjecture for spheres.
Also, the generalized lower bound conjecture (characterizing
the spheres for which $g_k=0$) is verified for
Bier spheres.
\end{itemize}
We are grateful to a referee for helpful remarks.

%%%%%%%%%%%%%%%%%%%%%%%%%%%%%%%%%%%%%%%%%%%%%%%%%%%%%%%%%%%%%%%%%%%%%

\section{Basic Definitions and Properties}

In this section we introduce our extension of Bier's construction to
bounded posets, and present some simple properties. We refer to
\cite{Stanley97} for background, notation and terminology relating to
posets and lattices.  Abstract simplicial complexes, order complexes, and
shellability are reviewed in \cite{Bjorner-topmeth}.  See \cite{Ziegler95}
for polytope theory.

All the posets we consider have finite length.  A poset is \emph{bounded}
if it has a unique minimal and maximal element; we usually denote these
by~$\hat0$ and~$\hat1$, respectively.  For $x\le y$, the \emph{length}
$\ell(x,y)$ is the length of of a longest chain in the \emph{interval}
$[x,y]=\{z\in P:x\le z\le y\}$.  A bounded poset is \emph{graded} if all
maximal chains have the same length.  A graded poset is \emph{Eulerian} if
every interval $[x,y]$ with $x<y$ has the same number of elements of odd
rank and of even rank.  An \emph{ideal} in $P$ is a subset $I\subseteq P$
such that $x\le y$ with $x\in P$ and $y\in I$ implies that $x\in I$. It is
\emph{proper} if neither $I=P$ nor $I=\emptyset$.  Our notation in the
following will be set up in such a way that all elements of $P$
named $x,x_i,x_i'$ are elements of the
ideal~$I\subset P$, while elements called $y,y_j,y_j'$ are 
in $P\sm I$.

\begin{defi}
  Let $P$ be a bounded poset of finite length and $I\subset P$ a proper
  ideal.  Then the poset $\bier(P,I)$ is obtained as follows: It
  consists of all intervals $[x,y]\subseteq P$ such that $x\in I$ and
  $y\notin I$, ordered by reversed inclusion, together with an additional
  top element $\hat 1$.
\end{defi}

Here \emph{reversed} inclusion says that $[x',y']\le[x,y]$
amounts to $x'\le x<y\le y'$.  The interval $I=[\hat0,\hat1]$ is
the unique minimal element of $\bier(P,I)$, so $\bier(P,I)$ is bounded.

One may observe that the construction of Bier posets has a curious formal
similarity to the $E_t$-construction of Paffenholz \& Ziegler as defined
in \cite{Z89}.  The study of posets of intervals in a given poset, ordered
by inclusion, goes back to a problem posed by Lindstr\"om
\cite{Lindstrom71}; see Bj\"orner
\cite{Bjorner-cwposets,Bjorner-antiprism} for results on interval posets
related to this problem.

\begin{lemma}\label{lemma:basics}
  Let $P$ be a poset and $I\subset P$ a proper ideal. 
\begin{itemize}\itemsep=-2pt
\item[\rm(i)] The posets $P$ and $\bier(P,I)$ have the same length~$n$.
\item[\rm(ii)] $\bier(P,I)$ is graded if
  and only if $P$ is graded.\\
  In that case, $\rank[x,y]=\rank_Px+(n-\rank_Py)$.
\item[\rm(iii)] The intervals of $\bier(P,I)$ are of the following two kinds:
\begin{eqnarray*}
[[x,y],\hat 1] &\cong& \bier([x,y],I\cap[x,y]) \\
{}[[x',y'],[x,y]]& =   & [x',x]\times[y,y']^{\mathrm{op}}, 
\end{eqnarray*}
where $[y,y']^{\mathrm{op}}$ 
denotes the interval $[y,y']$ with the opposite order.
\item[\rm(iv)] If $P$ is a lattice then $\bier(P,I)$ is a lattice.
\end{itemize}
\end{lemma}

\begin{proof}
  $\bier(P,I)$ is bounded. Thus for (iv) it suffices to show that meets
  exist in $\bier(P,I)$. These are given by $[x,y]\wedge[x',y']=[x\wedge
  x',y\vee y']$ and $[x,y]\wedge\hat 1=[x,y]$.  The other parts are
  immediate from the definitions.
\end{proof}

%%%%%%%%%%%%%%%%%%%%%%%%%%%%%%%%%%%%%%%%%%%%%%%%%%%%%%%%%%%%%%%%%%%%%

\section{Bier posets via stellar subdivisions}\label{sec:StellarSubdiv}

For any bounded poset $P$ we denote by $\ov P:=P\sm \{\hat0,\hat1\}$ the
proper part of $P$ and by~$\Delta(\ov P)$ the order complex of~$\ov P$,
that is, the abstract simplicial complex of all chains in $\ov P$ (see
\cite{Bjorner-topmeth}).

In this section we give a geometric interpretation of $\bier(P,I)$, by
specifying how its order complex may be derived from the order complex of
$P$ via stellar subdivisions.  For this, we need an explicit
description of stellar subdivisions for abstract simplicial complexes.
(See e.\,g.\ \cite[p.~15]{rourke72:_introd_piecew_linear_topol} for the
topological setting.)

\begin{defi}
  The \emph{stellar subdivision} $\stellar_F(\Delta)$ of a
  finite-dimensional simplicial complex $\Delta$ with respect to a
  non-empty face~$F$ is obtained by removing from $\Delta$ all faces that
  contain $F$ and adding new faces $G\cup\{v_F\}$ (with a new apex vertex
  $v_F$) for all faces $G$ that do not contain $F$, but such that $G\cup
  F$ is in the original complex.
\end{defi}

In the special case of a stellar subdivision of an edge~$E=\{v_1,v_2\}$,
this means that each face $G\in\Delta$ that contains $E$ is replaced by
three new faces, namely $(G\sm\{v_1\})\cup\{v_E\}$,
$(G\sm\{v_2\})\cup\{v_E\}$, and $(G\sm\{v_1,v_2\})\cup\{v_E\}$.  Note that
this replacement does not affect the Euler characteristic.

\begin{remark*}
  The stellar subdivisions in faces $F_1,\ldots,F_N$ commute, and
  thus may be performed in any order --- or simultaneously --- if and only
  if no two $F_i,F_j$ are contained in a common face $G$ of the complex,
  that is, if $F_i\cup F_j$ is not a face for~$i\neq j$.
\end{remark*}

\begin{theorem}\label{thm:ordercomplex}
  Let $P$ be a bounded poset of length $\ell (P)=n<\infty$, and let
  $I\subset P$ be a proper ideal.
  Then the order complex of~$\ov{\bier(P,I)}$ is obtained from the order
  complex of $\ov P$ by stellar subdivision on all edges of the form
  $\{x,y\}$, for $x\in\ov I$, $y\in\ov P\sm\ov I$, $x<y$. These stellar
  subdivisions of  edges $\{x,y\}$ must be performed in order 
of increasing length $\ell(x,y)$.
\end{theorem}

\begin{proof}  
  In the following, the elements denoted by $x_i$ or $x_i'$ will be
  vertices of $\ov P$ that are contained in $\ov I:=I\sm\{\hat0\}$,
  while elements denoted by $y_j$ or $y_j'$ are from $\ov P\sm\ov I$.  By
  $(x'_i,y_i')$ we will denote the new vertex created by subdivision of
  the edge $\{x_i',y_i'\}$.
  
  We have to verify that subdivision of all the edges of $\Delta(\ov P)$
  collected in the sets
  \[
  E_k\ :=\ \big\{\,\{x,y\} : x<y,\  
  \ell(x,y)=k,\ x\in\ov I,\ y\in \ov P\sm\ov I\,\big\}
  \]
  for $k=1,\dots,n-2$ (in this order) results in
  $\Delta(\ov{\bier(P,I)})$.  To prove this, we will explicitly describe
  the simplicial complexes $\Gamma_k$ that we obtain at intermediate
  stages, after subdivision of the edges in $E_1\cup\dots\cup E_k$. (The
  complexes $\Gamma_k$ are \emph{not} in general order complexes
  for~$0<k<n-2$.)%
\smallskip

\noindent{\itshape{\bfseries Claim.}
  After stellar subdivision of the edges of $\Delta(\ov P)$ in the edge
  sets $E_1,\dots,E_k$ (in this order), the resulting complex $\Gamma_k$
  has the faces
\begin{equation}\label{eq:gamma_k_simplices}
\{ x_1,x_2,\dots,x_r,
      (x'_1,y'_1),(x'_2,y'_2),\dots,(x'_t,y'_t),y_1,y_2,\dots,y_s\}
\end{equation}
where
\begin{compactenum}[\rm(i)~]
\item\label{eq:c1}~\vskip-20pt
\begin{equation*}
x_1<x_2<\dots<x_r< y_1<y_2<\dots< y_s\qquad(r,s\ge0)
\end{equation*}
must be a strict chain in $\ov P$ that may be empty,
but has to satisfy $\ell(x_r,y_1)\ge k+1$ if $r\ge1$ and $s\ge1$, while
\item\label{eq:c2}~\vskip-20pt
\begin{equation*}
[x'_t,y'_t]<\dots<[x'_2,y'_2]<[x'_1,y'_1]\qquad(t\ge0)
\end{equation*}
must be a strict chain in $\ov{\bier(P,I)}$ that may be empty,
but has to satisfy $\ell(x'_t,y'_t)\le k$ if $t\ge1$,
and finally
\item\label{eq:c3}~\vskip-20pt
\begin{equation*}
x_r\le x'_t\quad\text{ and }\quad y'_t\le y_1
\end{equation*}
must hold if both $r$ and $t$ are positive 
   resp.\ if both $s$ and $t$ are positive.
\end{compactenum}}
\smallskip

\noindent
The conditions \eqref{eq:c1}--\eqref{eq:c3} together imply that the chains
of $\Gamma_k$ are supported on (weak) chains in $\ov P$ of the form
\begin{equation*}%
\hat0<x_1<x_2<\dots<x_r\le x'_t\le\ldots\le x'_2\le x'_1<
                     y'_1\le y'_2\le\ldots\le y'_t\le y_1<y_2\ldots< y_s<\hat1.
\end{equation*}
In condition \eqref{eq:c3} not both inequalities can hold with equality,
because of the length requirements for \eqref{eq:c1} and~\eqref{eq:c2},
which for $r,s,t\ge1$ 
mandate that $\ell(x'_t,y'_t)\le k<\ell(x_r,y_1)$, and thus
$[x'_t,y'_t]\subset[x_r,y_1]$.  
\smallskip

We verify immediately that for $k=0$ the description of $\Gamma_0$ given
in the claim yields $\Gamma_0=\Delta(\ov P)$, since for $k=0$ the length
requirement for \eqref{eq:c2} does not admit any subdivision vertices.

For $k=n-2$ the simplices of $\Gamma_{n-2}$ as given by the claim cannot
contain both~$x_r$ and~$y_1$, that is, they all satisfy either $r=0$ or
$s=0$ or both, since otherwise we would get a contradiction between the
length requirement for \eqref{eq:c1} and the fact that any interval
$[x_r,y_1]\subseteq\ov P$ can have length at most $n-2$. Thus we obtain
that $\Gamma_{n-2}=\Delta(\ov{\bier(P,I)})$, if we identify the
subdivision vertices $(x'_i,y'_i)$ with the intervals $[x'_i,y'_i]$ in
$P$, the elements $x_i$ with the intervals $[x_i,\hat1]$ and the elements
$y_j\in\ov P\sm\ov I$ with the intervals~$[\hat0,y_j]$.

Finally, we prove the claim by verifying the induction step from $k$ to
$k+1$.  It follows from the description of the complex $\Gamma_k$ that 
no two edges in $E_{k+1}$ lie in the same facet.  Thus we can stellarly
subdivide the edges in $E_{k+1}$ in arbitrary order. Suppose the edge
$(x_r,y_1)$ of the simplex
\[
 \{x_1,\dots,x_{r-1},x_r,
  (x'_1,y'_1),(x'_2,y'_2),\dots,(x'_t,y'_t), y_1,y_2,\dots,y_s\}
\]
is contained in $E_{k+1}$. Then stellar subdivision yields the three new
simplices
  \begin{alignat*}{3}
    &\{ x_1,\dots,x_{r-1},&&
    (x_r,y_1),(x'_1,y'_1),(x'_2,y'_2),\dots,(x'_t,y'_t),
    &\,y_1,\,&y_2,\dots,y_s,\},
\\
    &\{ x_1,\dots,x_{r-1},&\,x_r,\,&
    (x_r,y_1),(x'_1,y'_1),(x'_2,y'_2),\dots,(x'_t,y'_t),
    &&    y_2,\dots,y_s,\},
\quad\text{and}\\
    &\{ x_1,\dots,x_{r-1},&&
    (x_r,y_1),(x'_1,y'_1),(x'_2,y'_2),\dots,(x'_t,y'_t),
    &&    y_2,\dots,y_s,\}.
  \end{alignat*}
  All three sets then are simplices of $\Gamma_{k+1}$, satisfying all the
  conditions specified in the claim (with $t$ replaced by $t+1$ and $r$ or
  $s$ or both reduced by~$1$).  Also all simplices of $\Gamma_{k+1}$ arise
  this way.  This completes the induction step.  
\end{proof}
  
  We can write down the subdivision map of the previous proof explicitly: The map
  \begin{align*}
    \pi:\|\Delta(\ov {\bier(P,I)})\|&\rightarrow\|\Delta(\ov P)\|\\%
    \intertext{is given on the vertices of $\Delta(\ov {\bier(P,I)})$ by} 
    [x,y]&\mapsto
    \begin{cases}
      \frac{1}{2}x+\frac{1}{2}y 
           \quad &\hat 0<x<y<\hat 1, x\in I, y\notin I\\
      x &\hat 0<x<y=\hat 1, x\in I, y\notin I\\
      y &\hat 0=x<y<\hat 1, x\in I, y\notin I
    \end{cases}\notag
  \end{align*}
  and extends linearly on the simplices of~$\Delta(\ov {\bier(P,I)})$.

\begin{cor}
$\|\Delta(\ov {\bier(P,I)})\|$ and  $\|\Delta(\ov P)\|$
are PL homeomorphic.\qed
\end{cor}

In the case where $P$ is the face poset of a regular PL CW-sphere or
-manifold, this implies that the barycentric subdivision of $\bier(P,I)$
may be derived from the barycentric subdivision of $P$ by stellar
subdivisions. In particular, in this case $\bier(P,I)$ is again the face
poset of a PL-sphere or manifold.

\begin{cor}\label{cor:sphere}
  If $P$ is the face lattice of a strongly regular PL CW-sphere then so
  is~$\bier(P,I)$.\qed
\end{cor}

\begin{cor}
  If $P$ is Cohen-Macaulay then so is~$\bier(P,I)$.
\end{cor}
\begin{proof}
This follows from the fact that Cohen-Macaulayness 
(with respect to arbitrary coefficients) is
a topological property \cite{munkres84:CM}.
\end{proof}

%%%%%%%%%%%%%%%%%%%%%%%%%%%%%%%%%%%%%%%%%%%%%%%%%%%%%%%%%%%%%%%%%%%%%%%%
\section{Eulerian Posets}\label{sec:Eulerian}

{}From now on we assume that $P$ is a graded poset of length $n$, and that
$I\subset P$ is a proper order ideal, with $\hat0_P\in I$ and
$\hat1_P\notin I$.  First we compute the $f$-vector
$f(\bier(P,I)):=(f_0,f_1,\dots,f_n)$, where $f_i$ denotes the elements of
\emph{rank} $i$ in the poset $\bier(P,I)$.  (This notation is off by $1$
from the usual convention in polytope theory, as in~\cite{Ziegler95}.)  By
definition we have $f_n(\bier(P,I))=1$ and
\begin{eqnarray*}
  f_i(\bier(P,I))& =& \#\{\,[x,y] :  x\in I, y\notin I,
  \rank_P x+n-\rank_P y=i\,\}
\end{eqnarray*}
for $0\le i\le n-1$. In particular, $f_0(\bier(P,I))=1$.

\begin{theorem}
  Let $P$ be an Eulerian poset and $ I\subset P$ a proper ideal. Then
  $\bier(P,I)$ is also an Eulerian poset.
\end{theorem}

\begin{proof}
  $\bier(P,I)$ is a graded poset of the same length as $P$ by Lemma
  \ref{lemma:basics}.
  Thus it remains to prove that all intervals of length $\ge1$ in
  $\bier(P,I)$ contain equally many odd and even rank elements.
  
  This can be done by induction. For length $\ell(P)\le1$ the claim is
  true. Proper intervals of the form $[[x,y],\hat 1]$ are, in view of
  Lemma \ref{lemma:basics}, Eulerian by induction.  Proper intervals of
  the form $[[x',y'],[x,y]]$ are Eulerian, since any product of
  Eulerian posets is Eulerian.  Finally the whole poset $\bier(P,I)$
  contains the same number of odd and even rank elements by the following
  computation:
  \begin{align}
    \sum_{i=0}^n (-1)^{n-i} f_i(\bier(P,I))
    &=1+\sum_{i=0}^{n-1} (-1)^{n-i} f_i(\bier(P,I))\notag\\
    &=1+\sum_{y\notin I}\sum_{{x\in I}\atop {x\le y}} 
               (-1)^{\rank(y)-\rank(x)}\notag\\
    &=1+\sum_{y\notin I}\sum_{x\le y}(-1)^{\rank(y)-\rank(x)}
             -\sum_{y\notin I}\sum_{{x\notin I}\atop {x\le y}}
             (-1)^{\rank(y)-\rank(x)}\label{1}\\
    &=1+0-\sum_{x\notin I}\sum_{x\le y} 
             (-1)^{\rank(y)-\rank(x)}\label{2}\\
    &=1+0-1=0\notag
  \end{align}
  where the first double sum in \eqref{1} is $0$ as $[\hat0_P,y]$ is
  Eulerian and $\rank(y)\ge 1$, and the double sum in \eqref{2} is~$-1$ as
  $[x,\hat1_P]$ is Eulerian and trivial only for $x=\hat1_P$.
\end{proof}

Alternatively, the result of the computation in this proof also follows
from the topological interpretation of $\bier(P,I)$ in the previous
section.

%%%%%%%%%%%%%%%%%%%%%%%%%%%%%%%%%%%%%%%%%%%%%%%%%%%%%%%%%%%%%%%%%%%%%

\section{Shellability of Bier spheres}

Now we specialize to Bier's original setting, where $P=B_n$ is the boolean
lattice of all subsets of the ground set $[1,n]=\{1,\dots,n\}$ (which may
be identified with the set of atoms of~$B_n$), ordered by inclusion.  We
will use notation like $[1,n]$ or $(x,n]$ freely to denote closed or
half-open sets of integers.
 
Any non-empty ideal in the boolean algebra $B_n$ can be interpreted as an
abstract simplicial complex with at most $n$ vertices, so we denote it
by~$\Delta$.

We get
\[
\bier(B_n,\Delta) \sm \{\hat 1\} \ \ =\ \ 
\{(B,C) :  \emptyset\subseteq B \subset C \subseteq
[1,n], B\in \Delta, C\notin \Delta\} 
\]
again ordered by reversed inclusion of intervals. We denote the facets of
$\bier(B_n,\Delta)$ by $(A;x):=(A,A\cup\{x\})\in\bier(B_n,\Delta)$ and the set
of all facets by $\mathcal F(\Delta)$.
 
The poset $\bier(B_n,\Delta)$ is the face lattice of a simplicial PL
$(n-2)$-sphere, by Corollary~\ref{cor:sphere}.  We will now prove a
strengthening of this, namely that $\bier(B_n,\Delta)$ is shellable. (As
is known, see e.\,g.\ \cite{Bjorner-topmeth}, shellability implies the
PL-sphericity for pseudomanifolds.)

\begin{theorem}\label{thm:shellable}
  For every proper ideal $\Delta\subset B_n$, the $(n-2)$-sphere
  $\bier(B_n,\Delta)$ is shellable.
\end{theorem}

\begin{proof}
  The shellability proof is in two steps. First we show that the rule
  \begin{align}\label{restriction}
    \begin{split}
    R:\mathcal F(\Delta)&\rightarrow \bier(B_n,\Delta)\\
    (A;x)&\mapsto(A\cap(x,n],A\cup[x,n]).
    \end{split}
  \end{align}
  defines a restriction operator in the sense of \cite{Bjoerner92}; that
  is, it induces a partition
  \[
  \bier(B_n,\Delta)\ \ =\ \ \biguplus_{(A;x)\in\mathcal
    F(\Delta)}[R(A;x),(A;x)],
  \]
  and the precedence relation forced by this restriction operator is
  acyclic. Thus, any linear extension of the precedence relation yields a
  shelling order.
  
  That the restriction operator indeed defines a partition can be seen as
  follows: Take any element $(B,C)\in\bier(B_n,\Delta)$. Set
  \begin{alignat*}{3}
  x :=&\min&\{y\in C\sm B :  B\cup(C&\,\cap&[1,y])\notin \Delta\}\\
      &\max&\{y\in C\sm B :  B\cup(C&\,\sm &[y,n])\in \Delta\} 
 \end{alignat*}
  and $A:=B\cup(C\cap[1,x))$. Then we have
\[
A\cap(x,n]\subseteq B \subseteq A \subset A\cup\{x\}\subseteq
  C\subseteq A\cup[x,n]
\]
  and thus $(B,C)$ is contained in $[R(A;x),(A;x)]$. 

To see that the intervals in the partition do not
  intersect we have to show that if
  both $R(A;x)\le(A';x')$ and $R(A';x')\le(A;x)$,
  then $(A;x)=(A';x')$. This is a special case of a more general
  fact we establish next, so  we do not give the argument here.

For any shelling order ``$\prec$'' that would induce $R$
as its
``unique minimal new face'' restriction operator
we are forced to require that 
if $R(A;x)\le (A';x')$ for two facets $(A;x)$ and $(A';x')$,
then $(A;x)\preceq(A';x')$.
By definition, $R(A;x)\le (A';x')$ means that
\begin{equation}\label{eq:prec1}
A\cap (x,n]\ \subseteq\ A'\ \subset\ A'\cup\{x'\}\ \subseteq A\cup[x,n],
\end{equation}
which may be reformulated as
\begin{equation}\label{eq:prec2}
(A \cup\{x \})_{>x}\ \subseteq\ A'\qquad\textrm{and}\qquad
(A'\cup\{x'\})_{<x}\ \subseteq\ A.
\end{equation}

We now \emph{define} the relation $(A;x)\prec(A';x')$
to hold if and only if \eqref{eq:prec2} holds
together with
\begin{equation}\label{eq:prec3}
(A \cup\{x \})_{\le x}\ \not\subseteq\ A'\qquad\textrm{and}\qquad
(A'\cup\{x'\})_{\ge x}\ \not\subseteq\ A.
\end{equation}
Note that our sets $A,A'$ belong to
an ideal which does not contain $A\cup\{x\},A'\cup\{x'\}$,
so \eqref{eq:prec3} applies if \eqref{eq:prec2} does. 

By the \emph{support} of $(A;x)$ we mean the set $A\cup\{x\}$.
The element $x$ of the support is called its \emph{root element}.

We interpret a relation $(A;x)\prec(A';x')$ as a \emph{step}
from $(A;x)$ to $(A';x')$.
The first conditions of \eqref{eq:prec2} and \eqref{eq:prec3}
say that
\begin{equation}\label{cond:1}
\begin{minipage}{130mm}
In each step, the elements that are deleted from the support
are $\le x$; moreover, we must either loose some element $\le x$
from the support, or we must choose $x'$ from $(A\cup\{x\})_{\le x}$,
or both.
\end{minipage}
\end{equation}
Similarly, the second conditions of~\eqref{eq:prec2} and~\eqref{eq:prec3}
say that 
\begin{equation}\label{cond:2}
\begin{minipage}{130mm}
In each step, the elements that are added to the support
are $> x$;\break moreover, we must either add some element $>x$
to the support, or we must keep $x$ in the support,
or both.
\end{minipage}
\end{equation}

Now we show that the transitive closure of the relation $\prec$
does not contain any cycles.
So, suppose that there is a cycle,
\[
(A_0;x_0)\ \prec\ (A_1;x_1)\ \prec\ \ldots\ \prec\ (A_k;x_k)=(A_0;x_0).
\]
First assume that not all root elements $x_i$ in this cycle are equal.
Then by cyclic permutation we may assume that $x_0$ is the smallest
root element that appears in the cycle, and that $x_1>x_0$. Thus $x_1$
is clearly not from $(A\cup\{x_0\})_{\le x_0}$,
so by Condition~\eqref{cond:1} we loose an element 
$\le x_0$ from the
support of $(A_0;x_0)$ in this step. But in all later steps the
elements we add to the support are $>x_i\ge x_0$, so the lost element
will never be retrieved. Hence we cannot have a cycle.

The second possibility is that all root elements in the cycle
are equal, that is, $x_0=x_1=\dots=x_k=x$.
Then by Conditions~\eqref{cond:1} and~\eqref{cond:2},
in the whole cycle we loose only elements $<x$ from the support,
and we add only elements $>x$. 
The only way this can happen is that, when we traverse the 
cycle,  no elements are lost and  none are added,  
so $A_0=A_1=\dots=A_k$. Consequently, there is no cycle.
\end{proof}

The relation defined on the set of all pairs $(A;x)$ with $A\subset[1,n]$
and $x\in[1,n]\sm A$ by \eqref{eq:prec2} alone does have cycles, such as
\[
(\{1,4\},2)\ \prec\ (\{1,4\},3)\ \prec\ (\{4\},1)\ \prec\ (\{1,4\},2).
\]
This is the reason why we also require condition~\eqref{eq:prec3}
in the definition of ``$\prec$''.
\smallskip

The shelling order implied by the proof of Theorem~\ref{thm:shellable}
may also be described in terms of a linear ordering.
For that we associate with each facet $(A;x)$ a vector
$\chi(A;x)\in\R^n$, defined as follows:
\[
\chi(A;x)_a\ :=\ 
\begin{cases}
-1           & \textrm{for } a\in (A\cup\{x\})_{\le x},\\
\hphantom{-}0& \textrm{for } a\notin\ A\cup\{x\},\qquad\textrm{and}\\
+1           & \textrm{for } a\in (A\cup\{x\})_{>x}.
\end{cases}
\]
With this assignment, we get that
$(A;x)\prec(A';x')$, as characterized by 
Conditions~\eqref{cond:1} and~\eqref{cond:2},
implies that $\chi(A;x)<_{\textrm{lex}}\chi(A';x')$.
Thus we have that lexicographic ordering on the $\chi$-vectors
induces a shelling order for every ``boolean Bier sphere.''

%%%%%%%%%%%%%%%%%%%%%%%%%%%%%%%%%%%%%%%%%%%%%%%%%%%%%%%%%%%%%%%%%%%%%

\section{\boldmath$g$-Vectors}

The $f$-vectors of triangulated spheres are of great
combinatorial interest. In this section we derive the basic relationship
between the $f$-vector of a Bier sphere $\bier(B_n,\Delta)$
and the $f$-vector of the underlying simplicial complex $\Delta$.
(Such an investigation had been begun in Bier's note \cite{Bier92}.)

In extension of the notation of Section~\ref{sec:Eulerian} let 
$f_i(\Delta)$ denote the number of sets of cardinality $i$
in a complex~$\Delta$. The \emph{$f$-vector} of a proper subcomplex
$\Delta\subset B_n$
is $f(\Delta)=(f_0,f_1,\dots,f_n)$, with $f_0=1$ and $f_n=0$.

Now let $\Gamma$ be a finite simplicial complex
that is pure of dimension $d=n-2$, that is, such that
all maximal faces have cardinality $n-1$.
(Below we will apply this to $\Gamma=\bier(B_n,\Delta)$.)
We define $h_i(\Gamma)$ by
\begin{equation}\label{def:h_i}
h_i(\Gamma)\ :=\ \sum_{j=0}^{n-1} (-1)^{i+j}\binom{n-1-j}{n-1-i}f_j(\Gamma)
\end{equation}
for $0\le i\le n-1$, and $h_i(\Gamma):=0$ outside this range.
Then, conversely
\[
f_i(\Gamma)\ =\ \sum_{j=0}^{n-1} \binom{n-1-j}{n-1-i}h_j(\Gamma).
\]
Finally, for $0\le i\le\lfloor \frac{n-1}2 \rfloor$
let $g_i(\Gamma):=h_i(\Gamma)-h_{i-1}(\Gamma)$,
with $g_0(\Gamma)=1$.

Now we consider the $f$-, $h$- and $g$-vectors of
the sphere $\Gamma=\bier(B_n,\Delta)$.
It is an $(n-2)$-dimensional shellable sphere on
$f_1(\Delta) + n - f_{n-1}(\Delta)$ vertices.
(So for the usual case of $f_1=n$ and $f_{n-1}=0$,
when $\Delta$ contains all the $1$-element subsets
but no $(n-1)$-element subset of~$[1,n]$, we
get a sphere on $2n$ vertices.)
In terms of the facets $(A;x) \in \mathcal F(\Delta)$
we have the following description of its $h$-vector:
\begin{equation}\label{restriction-h}
h_i (\bier(B_n,\Delta))\ =\ \#\{ (A;x) \in \mathcal F(\Delta):
 |\, A\cap(x,n]\,| + |\, [1,x) \sm A\,| =i\}
\end{equation}
for $0\le i\le n-1$.  This follows from the interpretation
of the $h$-vector of a shellable complex in terms of the restriction
operator as
\[
h_i(\bier(B_n,\Delta))\ =\ 
\#\{\,(A;x)\in \mathcal F(\Delta): \rank(R(A;x))=i\,\},
\]
see~\cite[p.~229]{Bjoerner92}, together with
equation~\eqref{restriction} and Lemma \ref{lemma:basics}\,(ii).

\begin{lemma}[Dehn-Sommerville equations]\label{DS} For $0\le i\le n-1$,
\[
h_{n-1-i} (\bier(B_n,\Delta))\ =\ h_i (\bier(B_n,\Delta)).
\]
\end{lemma}

\begin{proof}
It is a nontrivial fact that this relation is true
for \emph{any} triangulated $(n-2)$-sphere. 
However, in our situation it is a direct and elementary consequence of
equation \eqref{restriction-h}. 

Namely, neither the definition of the $h$-vector 
nor the construction of the Bier sphere depends
on the ordering of the ground set. Thus we can reverse
the order of the ground set $[1,n]$, to get that
\begin{equation}\label{restriction-h2}
h_i (\bier(B_n,\Delta))\ =\ \#\{ (A;x) \in \mathcal F(\Delta):
|\, A\cap[1,x) \,| + |\, (x,n]\sm A\,| =i\}.
\end{equation}
Thus a set $A$ contributes to $h_i(\bier(B_n,\Delta))$
according to \eqref{restriction-h} if and only if 
the complement of $A$ with respect to the $(n-1)$-element
set $[1,n]\sm\{x\}$ contributes to $h_{n-1-i}(\bier(B_n,\Delta))$
according to \eqref{restriction-h2}.
\end{proof}

The $g$-vector of $\bier(B_n,\Delta)$ has the following nice form.

\begin{theorem}\label{g-thm}
For all $i=0,\dots, \lfloor \frac{n-1}2 \rfloor$,
\[
g_i(\bier(B_n,\Delta))\ =\ f_i(\Delta)-f_{n-i}(\Delta).
\]
\end{theorem}

\begin{proof}
\newcommand{\caug}{\Delta^{\mathrm{aug}}}%
\newcommand{\bnnc}{\bier(B_{n+1},\caug)}%
Let $\caug$ be the same complex as $\Delta$, but viewed as sitting
inside the larger boolean lattice $B_{n+1}$.
We claim that
\begin{equation}\label{star}
h_i (\bnnc) \ =\  h_{i-1} (\bier(B_{n},\Delta)) + f_i(\Delta)
\end{equation}
for $0\le i\le n$. This is seen from equation \eqref{restriction-h2} as
follows.  The facets $(A;x)$ of $\bnnc$ that contribute to $h_i (\bnnc)$
are of two kinds: either $x\neq n+1$ or $x=n+1$.  There are $h_{i-1}
(\bier(B_n,\Delta))$ of the first kind and $f_i(\Delta)$ of the second.

Using both equation \eqref{star} and Lemma \ref{DS} twice we compute
\begin{eqnarray*}
g_i(\bier(B_n,\Delta))
&=& h_i      (\bier(B_n,\Delta))   - h_{i-1}(\bier(B_n,\Delta)) \\
&=& h_{n-1-i}(\bier(B_n,\Delta))   - h_{i-1}(\bier(B_n,\Delta)) \\
&=& h_{n-i}(\bnnc)-f_{n-i}(\Delta) - h_{i-1}(\bier(B_n,\Delta)) \\
&=& h_{i}  (\bnnc)-f_{n-i}(\Delta) - h_{i-1}(\bier(B_n,\Delta)) \\
&=& f_{i} (\Delta)-f_{n-i}(\Delta) .
\end{eqnarray*}\vskip-9mm
\end{proof}

\begin{cor}\label{depends}
  The face numbers $f_i(\bier(B_n,\Delta))$ of the Bier sphere depend only
  on $n$ and the differences $ f_i(\Delta) - f_{n-i}(\Delta)$.
\end{cor}

\begin{proof}
  The $g$-vector determines the $h$-vector (via Lemma \ref{DS}), which
  determines the $f$-vector.
\end{proof}

For example, if $n=4$ and $f(\Delta)=(1,3,0,0,0)$ or
$f(\Delta)=(1,4,3,1,0)$, then we get $g(\bier(B_4,\Delta))=(1,3)$ and
$f(\bier(B_4,\Delta))=(1,7,15,10)$.

\begin{theorem}
Every simplicial complex $\Delta\subseteq B_n$
has a subcomplex $\Delta'$ such that
\[
f_i(\Delta')=f_i(\Delta)-f_{n-i}(\Delta)
\]
for $0\leq i\leq\lfloor \tfrac{n}2\rfloor$ and
$f_i(\Delta')=0$ for $i>\lfloor \tfrac{n}2\rfloor$.
\end{theorem}

\begin{proof}
For any simplicial complex $\Delta$ in $B_n$, define the $d$-vector
by $d_i(\Delta)=f_i(\Delta)-f_{n-i}(\Delta)$
for $0\leq i\leq\lfloor \tfrac{n}2\rfloor$ and
$d_i(\Delta)=0$ for greater $i$. We shall find a subcomplex
$\Delta'\subseteq\Delta$ with $f_i(\Delta')=d_i(\Delta)$ for
all $i$.

Choose $\Delta'$ as a minimal subcomplex of $\Delta$ with the same
$d$-vector.
We must show that $f_i(\Delta')=0$ for all
$\lfloor \tfrac{n}2\rfloor<i\leq n$. Suppose that
there is a set $C\in\Delta'$ with $|C|>\tfrac{n}2$.
Then there is an involution $\pi:[1,n]\rightarrow[1,n]$,
i.\,e.\ a permutation of the ground set of order two, such that
\begin{equation}\label{eq:pi}
\pi(C)\supseteq [1,n]\sm C,
\end{equation}
where $\pi(C)$ is the image of $C$.
Now define $\varphi:B_n\rightarrow B_n$ by
$\varphi(B)=[1,n]\sm\pi(B)$ for all $B\subseteq[1,n]$.
Observe that $\varphi$ satisfies the following for all $B\subseteq[1,n]$:
\begin{itemize}
\item[(a)] $\varphi(\varphi(B))=B$,
\item[(b)] $B'\subseteq B\ \Rightarrow\ \varphi(B')\supseteq\varphi(B)$,
\item[(c)] $|B|+|\varphi(B)|=n$.
\end{itemize}
Let $K:=\{B\in\Delta'\;:\;\varphi(B)\in\Delta'\}$.  We claim that
$\Delta'\sm K$ is a simplicial complex with the same $d$-vector 
as~$\Delta'$.

First, we show that $\Delta'\sm K$ is a complex.  Let $B'\subseteq
B\in\Delta'\sm K$.  Then $B'\in\Delta'$ so we must show that
$B'\notin K$. Property (b) gives $\varphi(B')\supseteq\varphi(B)$, so we
get $B\notin K\Rightarrow\varphi(B)\notin\Delta'
\Rightarrow\varphi(B')\notin\Delta'\Rightarrow B'\notin K$.

Let $K_i=\{B\in K\,:\,|B|=i\}$ for $0\leq i\leq n$.  We have
$d_i(\Delta'\sm K)=(f_i(\Delta')-|K_i|)-(f_{n-i}(\Delta')-|K_{n-i}|)
=d_i(\Delta')-(|K_i|-|K_{n-i}|)$ for $0\leq i\leq\lfloor \tfrac{n}2\rfloor$.  We
must show that $|K_i|=|K_{n-i}|$ for all $i$. Property (a) gives that
$B\in K\Leftrightarrow\varphi(B)\in K$. Finally, property (c) gives that
$\varphi$ is a bijection between $K_i$ and $K_{n-i}$ for all $i$.

Fortunately, $K\neq\emptyset$ since
$\varphi(C)=[1,n]\sm\pi(C)\subseteq C$ by (\ref{eq:pi}), whence
$\varphi(C)\in\Delta'$ and $C\in K$. Thus we have found a strictly smaller
subcomplex of $\Delta'$ with the same $d$-vector --- a contradiction
against our choice of $\Delta'$.
\end{proof}

\begin{cor}\label{Kseq}
There is a subcomplex $\Delta'$ of $\Delta$ such that
\[
g_i(\bier(B_n,\Delta)) = f_i(\Delta')
\] 
for
$0\leq i\leq\lfloor\frac{n-1}2\rfloor$ and
$f_i(\Delta')=0$ for $i>\lfloor\frac{n-1}2\rfloor$.\qed
\end{cor}

It is a consequence of Corollary \ref{Kseq} that the $g$-vector $(g_0,g_1,
\dots, g_{\lfloor(n-1)/2\rfloor})$ of $\bier(B_n,\Delta)$ is a {\em
  $K$-sequence}, i.\,e., it satisfies the Kruskal-Katona theorem. This is of
interest in connection with the so called $g$-conjecture for spheres,
which suggests that $g$-vectors of spheres are $M$-sequences (satisfy
Macaulay's theorem).  $K$-sequences are a very special subclass of
$M$-sequences, thus $g$-vectors (and hence $f$-vectors) of Bier spheres
are quite special among those of general triangulated $(n-2)$-spheres on
$2n$ vertices. See \cite[Ch. 8]{Ziegler95} for details concerning $K$- and
$M$-sequences and $g$-vectors.

What has been shown also implies the following.
\begin{cor}\label{realize}
Every $K$-sequence $(1, n, \dots , f_k)$ with
$k\le \lfloor\frac{n-1}2\rfloor$ can be realized as the $g$-vector of a Bier sphere
with $2n$ vertices. \qed
\end{cor}

We need to review the definition of bistellar flips:
Let $\Gamma$ be a simplicial $d$-manifold. 
If $A$ is a\, $(d-i)$-dimensional face of\, $\Gamma$, $0\leq i\leq d$,
such that\, ${\rm link}_{\Gamma}(A)$ is the boundary\, $Bd(B)$ of an $i$-simplex $B$ 
that is not a face of $\Gamma$, then the operation $\Phi_A$ on $\Gamma$
defined by
\[
\Phi_A(\Gamma):=(\Gamma \backslash (A*Bd(B)))\cup (Bd(A)*B)
\]
is called a\, \emph{bistellar $i$-flip}.
Then $\Phi_A(\Gamma)$ is itself a simplicial $d$-manifold, homeomorphic
to $\Gamma$, and if
$0\le i\le \lfloor\tfrac{d-1}2 \rfloor$, then
\begin{align}\label{g-flip}
\begin{split}
g_{i+1}(\Phi_A(\Gamma)) & =  g_{i+1}(\Gamma) +1\\
g_{j}(\Phi_A(\Gamma)) & =  g_{j}(\Gamma)\, \quad\mbox{for all\,\,\,\, $j\neq i+1.$}
\end{split}
\end{align}
Furthermore, if\, $d$ is even and\, $i=\frac{d}{2}$, 
then\, $g_{j}(\Phi_A(\Gamma)) = g_{j}(\Gamma)$\, for all\, $j$.
See Pachner~\cite[p.~83]{Pachner}.
\smallskip

It follows from Corollary \ref{Kseq} that 
$g_k(\bier(B_n,\Delta))\ge 0$. The case of equality is characterized as follows.

\begin{cor}\label{LBC}
For $2\le k\le\lfloor \frac{n-1}2 \rfloor$, the following are equivalent:
\begin{enumerate}
\item[\rm (1)]$ g_k(\bier(B_n,\Delta))\ =0,$
\item[\rm (2)]  $f_k(\Delta)=0$ or $f_{n-k}(\Delta)=\binom{n}{i}$,
\item[\rm (3)] $\bier(B_n,\Delta)$ is obtained from the boundary
complex of the $(n-1)$-simplex via a sequence of bistellar $i$-flips,
with $i\le k-2$ at every flip.
\end{enumerate}
\end{cor}

\begin{proof}
$(1)\Rightarrow (2):$
Consider the bipartite graph
\newcommand{\Gnk}{G_{n,k}}
$\Gnk$ whose edges are the pairs $(A,B)$ such that $A$ is a $k$-element
subset, $B$ is an $(n-k)$-element subset of $[1,n]$, and $A\subset B$,
where the inclusion is strict since $k<n-k$. Then $\Gnk$ is a regular
bipartite graph (all vertices have the same degree), so by standard
matching theory $\Gnk$ has a complete matching.  The restriction of such a
matching to the sets~$B$ in~$\Delta$ gives an injective mapping
$\Delta_{n-k} \rightarrow \Delta_{k}$ from $\Delta$'s faces of cardinality
$n-k$ to those of cardinality $k$.

Equality $f_{n-k}(\Delta) = f_{k}(\Delta)$ implies that $\Gnk$ consists of
two connected components, one of which is induced on $\Delta_{n-k} \cup
\Delta_{k}$.  A nontrivial such splitting cannot happen since $\Gnk$ is
connected, so either $\Delta_{n-k}$ and $\Delta_{k}$ are both empty, or
they are both the full families of cardinality $\binom{n}{k}$.
\smallskip

$(2)\Rightarrow (3):$
As shown in  \cite{Bier92} and
\cite[Sect.~5.6]{Matousek03}, adding an $i$-dimensional face to
$\Delta$ produces a bistellar $i$-flip in $\bier(B_n,\Delta)$.
Now, $\Delta$ can be obtained from the empty complex by adding
$i$-dimensional faces, and here all $i\le k-2$ if
 $f_k(\Delta)=0$ (meaning that there are no faces of
dimension $k-1$ in $\Delta$). 
The case when $f_{n-k}(\Delta)=\binom{n}{i}$
is the same by symmetry.
\smallskip

$(3)\Rightarrow (1):$
This follows directly from (\ref{g-flip}), since the boundary 
of the $(n-1)$-simplex
has $g$-vector $(1,0,\dots,0)$.
\end{proof}

A convex polytope whose boundary complex is obtained from the boundary
complex of the $(n-1)$-simplex via a sequence of bistellar $i$-flips,
with $i\le k-2$ at every flip, is called \emph{$k$-stacked}. The {\em
  generalized lower bound conjecture} for polytopes maintains that
$g_k=0$ for a polytope if and only if it is $k$-stacked. This is still
open for general polytopes. See McMullen \cite{McMullen02} for a
recent discussion. Corollary \ref{LBC} shows that it is valid for
those polytopes that arise via the Bier sphere construction.

%%%%%%%%%%%%%%%%%%%%%%%%%%%%%%%%%%%%%%%%%%%%%%%%%%%%%%%%%%%%%%%%%%%%%

\section{Further Observations}

\subsubsection*{6.1 Many spheres}

In the introduction we remarked that the (isomorphism classes) of Bier
spheres are numerous, in fact so numerous that one concludes that most
of them lack convex realization. To show this, it suffices to consider
Bier spheres $\bier(B_n,\Delta)$ for complexes $\Delta$ that contain 
all sets $A\subset[1,n]$ of size $|A|\le \lfloor \frac{n-1}2\rfloor$,
a subcollection of the sets of size 
$|A|= \lfloor\frac{n-1}2\rfloor+1=\lfloor\frac{n+1}2\rfloor$, 
and no larger faces.
Equivalently, $\Delta$ is a complex of dimension at most 
$ \lfloor\frac{n-1}2\rfloor$ with complete 
$(\lfloor\frac{n-1}2\rfloor-1)$-skeleton.
There are $\binom{n}{\lfloor (n+1)/2\rfloor}=\binom{n}{\lfloor n/2 \rfloor}$ 
elements in the
$\lfloor\frac{n+1}2\rfloor$-level of $B_n$; thus there are at least
\[
\frac{2^{\binom{n}{\lfloor n/2\rfloor}}}{(2n)!}
\ \ \sim\ \ 
\frac{2^{2^n/\sqrt{n}}}
{(\frac{2n}{e})^{2n}}
\]
combinatorially non-isomorphic such Bier spheres (where our rough
approximation ignores polynomial factors). On the other hand
there are at most $2^{8n^3+O(n^2)}$ combinatorially non-isomorphic
simplicial polytopes on $2n$ vertices (see Goodman \& Pollack \cite{GP86},
Alon \cite[Thm.~5.1]{Alon86}).

It is interesting to contrast this with all the ways in which these
``numerous'' spheres are very special: They are shellable,
their $g$-vectors are $K$-sequences,
and for even $n$ we even get numerous
``nearly neighborly'' examples (as discussed below).
Another construction of ``numerous'' shellable spheres is known from the
work of Kalai \cite{Kalai88} and Lee \cite{Lee00}.

Though we have defined the construction of a Bier poset for 
arbitrary posets and have shown that the construction produces sphere
lattices from sphere lattices, it remains an open problem how to
extend the Bier construction to obtain numerous 
simplicial/shellable
$(n-2)$-spheres with more than $2n$ vertices.

\subsubsection*{6.2 Centrally symmetric and nearly neighborly spheres}

Let $\Gamma$ be a triangulated $(n-2)$-sphere on $2m$ vertices.
The sphere $\Gamma$ is \emph{centrally symmetric} if it has a 
symmetry of order two which fixes no face; that is,
if there is a fixed-point-free involution on its set $V$ of vertices
such that (i) for every face $A$ of~$\Gamma$ also $\alpha(A)$ is a 
face, and (ii) $\{x, \alpha (x)\}$ is not a face, for all $x\in V$.
A subset $A\subseteq V $ is \emph{antipode-free} if it
contains no pair $\{x, \alpha (x)\}$, for $x\in V$.  

A centrally symmetric sphere $\Gamma$ with involution $\alpha$ 
is \emph{$k$-nearly neighborly} 
if all antipode-free sets $A\subseteq V$ of size 
$|A|\le k$ are faces of~$\Gamma$.
Equivalently, $\Gamma$ must contain the $(k-1)$-skeleton of the
$m$-dimensional hyperoctahedron (cross-polytope). 
$\Gamma$ is \emph{nearly neighborly} if it is 
$\lfloor\frac{n-1}2\rfloor$-nearly neighborly.

Thus $k$-nearly neighborliness is defined only for centrally-symmetric
spheres. In the case $k\ge2$ the involution $\alpha$ is 
uniquely determined by the condition $\{x, \alpha (x)\}\notin\Gamma$.

The concept of nearly neighborliness 
for centrally symmetric spheres has been studied for centrally
symmetric $(n-1)$-polytopes, where $\alpha$ is of course the map $x\mapsto -x$.
For instance, work of Gr\"unbaum, McMullen and Shephard, Schneider,
and Burton shows that there are severe restrictions to $k$-nearly
neighborliness in the centrally symmetric polytope case, while
existence of interesting classes of nearly neighborly spheres was
proved by Gr\"unbaum, Jockusch, and Lutz;
see~\cite[p.~279]{Ziegler95} and \cite[Chap.~4]{Lutz-diss}.  

Nearly neighborly Bier spheres arise as follows.
(In the following, only the special case $m=n$,
of an $(n-2)$-sphere with $2n$ vertices, will occur.)
\begin{prop}\label{prop:cs}
  If $A\in\Delta\Longleftrightarrow[1,n]\sm A\notin\Delta$, then
  $\bier(B_n,\Delta)$ is centrally symmetric.
\end{prop}
\begin{proof}
  The involution $\alpha$ is given by the pairing 
  $[\{x\}, \hat{1}] \longleftrightarrow [\hat{0}, [1,n]\sm \{x\}].$
\end{proof}
\begin{prop}
Let $1< k \le \lfloor\frac{n-1}2\rfloor$.
The Bier sphere $\bier(B_n,\Delta)$ is a $k$-nearly neighborly 
$(n-2)$-sphere with $2n$ vertices if and only if 
\begin{enumerate}
\item[\rm(i)]
$A\in\Delta\Longleftrightarrow[1,n]\sm A\notin\Delta$,
for all $A\subseteq[1,n]$,
\item[\rm(ii)]
$B\in\Delta$, for all $B\subseteq[1,n]$, $|B|\le k$ \\
(and thus $C\notin\Delta$ for all $C\subseteq[1,n]$, $|C|\ge n-k$).
\end{enumerate}
\end{prop}
\begin{proof}
The Bier sphere $\bier(B_n,\Delta)$ has $2n$ vertices if and only if 
$\Delta\subset2^{[1,n]}$ is a complex that contains all
subsets of cardinality $1$ and no subsets of cardinality $n-1$.
The antipode-free vertex sets of cardinality $k$ then correspond
to intervals $[B,C]\subseteq B_n$ 
such that $|B|+(n-|C|)=k$.
A set $B$ is the minimal element of such an interval if and only if $|B|\le k$,
while $C$ is a  maximal element for $|C|\ge n-k$.
\end{proof}  

Combining these two propositions we obtain a large number of 
even-dimensional nearly
neighborly centrally symmetric Bier spheres.
Indeed, in the case of even $n$ we get at least
\[
{\frac{{2^{\frac12\binom{n}{\lfloor n/2\rfloor}}}}{(2n)!}}
\]
non-isomorphic spheres, from the complexes $\Delta$ which contain
all sets of size $A<\frac n2$,
and exactly one set from each pair of sets $A$ and $[1,n]\sm A$
of size $|A|=\frac n2$.

On the other hand, for odd $n$ (that is, in the case of an 
odd-dimensional sphere, or an even-dimensional polytope,
where the ``nearly neighborliness condition'' is stronger
and hence more interesting) only one instance of a 
nearly neighborly centrally symmetric Bier $(n-2)$-sphere 
with $2n$ vertices is obtained; namely, 
for $\Delta=\{A\subset[1,n]:|A|\le\lfloor\frac n2\rfloor\}$.

\end{document}